\documentclass[11pt]{article}

\usepackage{latexsym}
\usepackage{amsmath}
\usepackage{amssymb}
\usepackage{amscd}
\usepackage{amsthm}

\theoremstyle{plain}
 \newtheorem{theorem}{Theorem}[section]
 
 \newtheorem{lemma}[theorem]{Lemma}

\theoremstyle{definition}
 \newtheorem{definition}[theorem]{Definition}

\theoremstyle{remark}
 \newtheorem{remark}[theorem]{Remark}

\renewcommand{\qedsymbol}{$\blacksquare$}

\newcommand{\R}{\mathbb{R}}

\DeclareMathOperator{\projdim}{projdim}

\DeclareMathOperator{\depth}{depth}

\begin{document}

\title{L. Szpiro's conjecture on Gorenstein algebras in codimension 2} 

\author{Christian B\"ohning}

\date{ 2003}

\maketitle

\begin{abstract}
A Gorenstein $A-$algebra $R$ of codimension $2$ is a perfect finite $A-$algebra such that $R\cong \mathrm{Ext}^2_A(R ,A)$ holds as $R-$modules, $A$ being a Cohen-Macaulay local ring with $\mathrm{dim}\, A -\mathrm{dim}_A \, R=2$.\\ I prove a structure theorem for these algebras improving on an old theorem of M. Grassi [Gra]. Special attention is paid to the question how the ring structure of $R$ is encoded in its Hilbert resolution. It is shown that $R$ is automatically a ring once one imposes a very weak depth condition on a determinantal ideal derived from a presentation matrix of $R$ over $A$. Furthermore, the interplay of Gorenstein algebras and Koszul modules as introduced by M. Grassi is clarified. I include graded analogues of the afore-mentioned results when possible.\\
Questions of applicability to the theory of surfaces of general type (namely, canonical surfaces in $\mathbb{P}^4$) have served as a guideline in these commutative algebra investigations.
\end{abstract}

\setcounter{section}{-1}

\section{Introduction and statement of results}
Motivation for investigating the types of questions treated in the present article sprang from two different, though closely related sources, the first algebro-geometric, the second one purely algebraic in spirit. Though in this paper I will restrict myself to touching upon the latter only, both themes can, I think, be best understood in conjunction, so I will briefly discuss them jointly in the introduction.\\
From the point of view of algebraic geometry, the perhaps earliest traces of the story may be located in the book [En] by F. Enriques where he treats the structure of canonical surfaces in $\mathbb{P}^4$ with $q=0,\: p_g=5 ,\: K^2=8$ and $9$ (cf. loc. cit., p. 284ff. ; they are the complete intersections of type $(2,4)$ and $(3,3)$; the word ``canonical'' means that for those surfaces the $1-$canonical map gives a birational morphism onto the image in $\mathbb{P}^4$). Subsequently the case $K^2=10$ was solved by C. Ciliberto (cf. [Cil]) using liaison arguments, and along different lines, D. Ro\ss berg (cf. [Ro\ss ]) tackled the problem for $K^2=11$ and $=12$, the question being open for higher values of $K^2$.\\
Translated into the language of commutative and homological algebra, the first difficulty to face in an attempt to treat cases with higher $K^2$ is to find a satisfactory structure theorem for Gorenstein algebras in codimension $2$; roughly, these are finite $A-$algebras $R$ ($A$ some ``nice'' base ring) with $R\cong \mathrm{Ext}^2_A(R, A)$, possibly up to twist if the base ring is graded (cf. section $2$ below for precise definitions). The connection with the above surfaces is established by remarking that their canonical rings are codimension $2$ Gorenstein algebras over the homogeneous coordinate ring of $\mathbb{P}^4$.\\
With regard to a structure theorem, the above qualifier ``satisfactory'' means precisely that one should be able to tell from practically verifiable and non-tautological conditions how the Hilbert resolution of $R$ over $A$ encodes (1) the ``duality'' $R\cong \mathrm{Ext}_A^2(R,A)$ and (2) the fact that $R$ has not only an $A-$module structure, but also a ring structure. Whereas (1) is by now fairly well understood, (2) is not, and the main purpose of this paper is to show how (2) can be disposed with. Let me mention at this point that, in the course of his investigations concerning low rank vector bundles on projective spaces, Lucien Szpiro was presumably the first to point out the need for and formulate some conjectural statements concerning a good structure theorem for Gorenstein algebras in codimension $2$, which is why his name appears in the title of this paper.\\
Since then, quite a good deal of work has been done on this problem, in geometric and algebraic guises. Let me therefore give some perspective on its history: In [Cat2] canonical surfaces in $\mathbb{P}^3$ are studied (from a moduli point of view) via a structure theorem proved therein for Gorenstein algebras in codimension 1. It is shown that the duality $\R\cong \mathrm{Ext}_A^1(R , A)$ for these algebras translates into the fact that the Hilbert resolution of $R$ can be chosen to be self-dual; moreover, that the presence of a ring structure on $R$ is equivalent to a (closed) condition on the Fitting ideals of a presentation matrix of $R$ as $A-$module (the so-called ``ring condition'' or ``rank condition'' or ``condition of Rouch\'e-Capelli'', abbreviated R.C. in any case). These ideas were developed further and generalized in [M-P] and [dJ-vS] (within the codimension 1 setting). In particular the latter papers show that R.C. can be rephrased in terms of annihilators of elements of $R$ and gives a good structure theorem also in the non-Gorenstein case. On the base of all this, M. Grassi turned attention to the codimension $2$ setting. In [Gra] he isolated the abstract kernel of the problem and proved that also for codimension $2$ Gorenstein algebras the duality $R\cong \mathrm{Ext}^2_A(R ,A)$ is equivalent to $R$ having a self-dual resolution. He introduced the concept of Koszul modules which provide a nice framework for dealing with Gorenstein algebras and also proposed a structure theorem for the codimension $2$ case. Unfortunately, as for the question how the ring structure of a codimension $2$ Gorenstein algebra is encoded in its Hilbert resolution, the conditions he gives are tautological and (therefore) too complicated (although they are necessary and sufficient). More recently, D. Eisenbud and B. Ulrich (cf. [E-U]) re-examined the ring condition and gave a generalization of it which appears to be more natural than the direction in which [Gra] is pointing. But essentially, they only give sufficient conditions for $R$ to be a ring, and these are not fulfilled in the applications to canonical surfaces one has in mind. More information on the development sketched here can be found in [Cat4]. For a deeper study of that part of the story that originates from the duality $R\cong \mathrm{Ext}^2_A(R ,A)$ and its effects on the symmetry properties of the Hilbert resolution of $R$, as well as for a generalization of this to the bundle case cf. [E-P-W].\\
In the present article, on the contrary, the approach to the problem of detecting the ring structure of $R$ in its Hilbert resolution is based on a philosophy already present in [Cat4], cf. p. 48: The ring condition R.C. is automatical under some mild extra condition which, moreover, works well in the applications to canonical surfaces (the condition corresponds to the requirement that the canonical image in $\mathbb{P}^4$ should have only isolated singularities). This is the content of theorem 1.1 (and its graded analogue theorem 1.3) below. The idea of proof is marvellously simple: One considers the sheaf $\mathcal{R}$ associated to $R$ on $X:=\mathrm{Spec}\, A$ and uses the fact that ``the locus where $\mathcal{R}$ is not known to be a sheaf of rings a priori is small'', i.e. if $Y$ is the support of $\mathcal{R}$, there is an open $U\subset Y$ such that the complement $Z:=Y\backslash U$ has codimension $2$ in $Y$ and $\mathcal{R}|_U$ is known to be a sheaf of rings from elementary considerations. Then using the fact that $R$ is Cohen-Macaulay and some local cohomology computations, one can easily check that the ring structure extends from $\mathcal{R}|_U$ to $\mathcal{R}$. The situation is thus very much reminiscent of that encountered in the familiar Hartogs theorem in several complex variables.\\
Theorem 1.1 works without the Gorenstein condition on $R$ which will enter only in section 2. There I basically only combine the results from section 1 with the symmetry statements known for the resolutions of Gorenstein algebras from [Gra]. I have chosen to dwell on the proofs of the latter, partly because occasionally minor simplifications could be made, partly for the sake of completeness. Theorem 2.4 (resp. theorem 2.5) contains the characterization of Gorenstein algebras in codimension 2.\\
Section 3 is included to provide a strengthening of theorem 2.4 (the local case of the structure theorem). It is also meant to clarify the relationship between Gorenstein algebras in codimension 2 and Koszul modules as introduced by M. Grassi in [Gra]. It is shown that a Gorenstein algebra over a local ring (always assumed to be Cohen-Macaulay) admits a Gorenstein-symmetric resolution which is at the same time of Koszul module type (see section 3 for precise statements).\\
Finally, there is reasonable hope that the techniques developed in this paper may facilitate the study of canonical surfaces in $\mathbb{P}^4$ with $p_g=5,\: q=0 ,\: K^2\ge 13$. In fact I propose to apply them as well to questions of existence as to the problem of describing their moduli. Tangible results have already been obtained and will be published elsewhere.\\
It is a pleasure for me to thank Fabrizio Catanese for posing the problem and many useful discussions.    

\section{How the ring structure is encoded}
The following result should be viewed as a rather general extension theorem giving
conditions under which a module that is a ring "in codimension 1" is already itself a
ring, and giving a description of the resulting structure. 
\begin{theorem}
Let $A$ be a Cohen-Macaulay local ring with maximal ideal $\mathfrak{m}$ and residue
field $k$ and let $R$ be a finite $A-$module.\\
Let a length 2 minimal free resolution
\[
\begin{CD}
0 @>>> F_2 @>{\psi }>> F_1 @>{\varphi }>> F_0 @>p>> R @>>> 0
\end{CD}
\]
of $R$ be given. View $\varphi$ as a presentation matrix of $R$ and let $\varphi '$
be the matrix $\varphi$ with first row erased. Denote by $I$ resp. $I'$ the zeroth
Fitting ideals (=ideals of maximal minors) of $\varphi$ resp. $\varphi'$. \\
Suppose that $\depth (\mathrm{Ann}_A\, R , A)(= \mathrm{codim}_{A} R
)=2$, the maximum possible in view of the inequality $\depth (\mathrm{Ann}_A\, R ,
A)\le \projdim_A R$.\\
Then, if moreover the following condition holds
\begin{itemize}
\item[($\heartsuit$)] $\depth \, I' \ge 4$,
\end{itemize}  
it follows that $R\cong \mathrm{Hom}_{A_Y} (C ,C )$, where $A_Y:=A/\mathrm{Ann}_A\, R$
and $C:=\mathrm{Hom}_{A_Y} (R, A_Y)$ is the so-called conductor of $R$ into $A_Y$; and
$R$ is given a (commutative) ring structure if we define multiplication as composition
of endomorphisms of $C$. The identity element of $R$ is the image under $p$ of
the generator of $F_0$ corresponding to the first row of $\varphi$.  
\end{theorem}

\begin{remark}
If $R$ has a symmetric minimal free resolution of length 2, i.e. a resolution of the
form
\begin{displaymath}
\begin{CD}
0@>>> (A^n)^{\vee} @>{\psi={{-\beta^{\vee}}\choose {\alpha^{\vee}}}}>> A^n\oplus
(A^n)^{\vee} @>{ \varphi=(\alpha\, \beta )}>> A^n @>p>> R @>>> 0 ,
\end{CD}
\end{displaymath}
then the condition $\depth ( \mathrm{Ann}_A\, R , A) =2$ follows from the
Eisenbud-Buchsbaum acyclicity criterion (cf. [Ei], thm. 20.9); for then the zeroth Fitting ideals of $\psi$
and $\varphi$ agree and we must have $\depth \, I\ge 2$ for this complex to be exact.
But $\mathrm{rad}\, I =\mathrm{rad}\, \mathrm{Ann}_A \, R$ whence, in view of $\depth (
\mathrm{Ann}_A\, R , A) \le \projdim_A R$, $\depth ( \mathrm{Ann}_A\, R , A) =2$.
\end{remark}

\begin{proof}
The assumption that $A$ be Cohen-Macaulay implies that $R$ is a Cohen-Macaulay
$A-$module, too; the argument is that the Auslander-Buchsbaum formula $\projdim \, R+
\depth\, R=\depth \, A=\dim A$ (the latter because $A$ is Cohen-Macaulay) gives
$\depth \, R= \dim A -2$, whereas $2=\depth \,(\mathrm{Ann}_A\, R ,A)= \dim A -\dim_A
\, R$ (where again enters the hypothesis that $A$ is Cohen-Macaulay whence the depth of
an ideal is given by its codimension), and thus $\depth\, R=\dim_A R$, i.e. $R$ is
Cohen-Macaulay. \\
It will be convenient to introduce the following notation; let
\begin{quote}
\begin{tabular}{|l}
$X:=\mathrm{Spec}\, A$,\\
$Y$:= the codimension 2 closed subscheme of $X$ associated to $\mathrm{Ann}_A\, R$,\\
$Z\subset Y\subset X$ the codimension 2 (in $Y$) closed subscheme defined by $I'$,\\
$U:=Y-Z$ the open complement, $j: U\hookrightarrow Y$ the inclusion,\\
$\mathcal{R}$ resp. $\mathcal{C}$ the sheaves associated to $R$ resp. $C$ on $Y$.
\end{tabular}
\end{quote}
Now $\mathcal{R}|_U\cong \mathcal{O}_Y|_U$ (since the localization of the presentation matrix $\varphi'$ for $R/A_Y$ has invertible maximal minors), which is what I paraphrased in the
sentence preceding theorem 1.1 saying that $R$ is a ring "in codimension 1". One has
the exact sequence relating local and global cohomologies (cf. [Groth], prop. 2.2)
\[
0\longrightarrow H^0_Z(Y, \mathcal{R}) \longrightarrow H^0(Y,\mathcal{R})
\longrightarrow H^0(U,
\mathcal{R})
\longrightarrow H^1_Z(Y,\mathcal{R})\longrightarrow 0 ;
\]
On the other hand, since $R$ is Cohen-Macaulay, $\depth (I', R)= \dim R -\dim R/I'R\ge
2$, whence (cf. [Groth], thm. 3.8) $H^0_Z(Y, \mathcal{R})=H^1_Z(Y, \mathcal{R})=0$ and
thus $R\cong \Gamma (U, \mathcal{O}_Y|_U)$ (as $A_Y-$modules). Thus, via this
isomorphism, we already know that $R$ carries a structure of $A_Y-$(or if you like
$A-$)algebra.

To prove that $R\cong \mathrm{Hom}_{A_Y}(C,C)$ I first have to digress on some
general points concerning the structure of $R$ and $C$ (cf. [E-U]): Let $e\in R$
denote the image under $p$ of the $A-$module generator of $F_0$ corresponding to the
first row of $\varphi$; then $\varphi'$ is a presentation matrix (over $A$) of
$R/A_Ye$ whence $I'\subset \mathrm{Ann}_A\, R/A_Ye$ or $(I'\cdot A_Y)R\subset A_Ye$.
But we have already seen that $\depth (I', R)\ge 2$ and thus there is an element
$d\in (I'\cdot A_Y)\subset A_Y$ which is a nonzerodivisor on $R$ (therefore also on
$A_Y$) with $d R\subset A_Y e \subset R$. That is, $R$ is what is called a
\emph{finite birational} $A-$module in [E-U]. Incidentally, this implies that the
algebra structure on $R$ is unique since it is a subalgebra of $R[d^{-1}]=A[d^{-1}]$.
This being said, the conductor $\mathrm{Hom}_{A_Y}(R, A_Y)$ naturally identifies with
the ideal $\mathrm{Ann}_{A_Y}(R/A_Y e)$ in $A_Y$ via the prescription that $a\in
\mathrm{Ann}_{A_Y}(R/A_Y e)$ should be identified with multiplication by $a$ on $R$ (
to get a map back, send $c\in\mathrm{Hom}_{A_Y}(R ,A_Y)$ to $c(d)/d$). In fact, the
conductor is also an ideal in $R$.\\
This last remark allows one to immediately conclude that
$R^{\ast\ast }\equiv\mathrm{Hom}_{A_Y}(C,A_Y)=\mathrm{Hom}_{A_Y} (C,C)$, cf. [Cat4],
lemma 5.3 (here $\ast$ denotes the dualizing functor with respect to $A_Y$): Indeed,
$\mathrm{Hom}_{A_Y}(C,C)\subset \mathrm{Hom}_{A_Y}(C,A_Y)$ being clear, let $\xi \in
\mathrm{Hom}_{A_Y}(C ,A_Y)$ be given; let $c$ be in $C$ arbitrary. We have to show
that $\xi c \in  C$. But if $r$ is in $R$, then $c r$ is in $C$, since $C$ is an
ideal in $R$, thus  $\xi c r$ is in $A_Y$ and hence $\xi c \in C$.\\
We are now ready to conclude that $R\cong \mathrm{Hom}_{A_Y}(C,C)$. Indeed, we know
already that $\mathcal{R}|_U \cong \mathcal{H}om_{\mathcal{O}_Y} (\mathcal{C} ,
\mathcal{C} )|_U$, since also $\mathcal{C}|_U\cong \mathcal{O}_Y|_U$. On the other
hand, we know $\mathcal{R}\cong j_{\ast }j^{\ast} \mathcal{R}$ (from the exact
sequence 
\[
0 \longrightarrow \mathcal{H}^0_Z(\mathcal{R}) \longrightarrow \mathcal{R}
\longrightarrow j_{\ast }j^{\ast} \mathcal{R} \longrightarrow
\mathcal{H}^1_Z(\mathcal{R})\longrightarrow 0
\] 
(cf. [Groth], cor. 1.9) and $\depth (I' ,R)\ge 2$). From the short exact
sequence
\[
0\longrightarrow \mathcal{C} \longrightarrow \mathcal{O}_Y \longrightarrow
\mathcal{O}_Y/\mathcal{C} \longrightarrow 0
\]
we get the long exact sequence
\[
0\longrightarrow \mathcal{H}^0_Z(\mathcal{C}) \longrightarrow
\mathcal{H}^0_Z(\mathcal{O}_Y) \longrightarrow
\mathcal{H}^0_Z(\mathcal{O}_Y/\mathcal{C})\longrightarrow
\mathcal{H}^1_Z(\mathcal{C})\longrightarrow
\ldots 
\]
and $\mathcal{H}^0_Z(\mathcal{O}_Y)=0$ because e.g. $d$ above is a nonzerodivisor on
$A_Y$, that is $\depth (I' , A_Y)\ge 1$. Thus also $\mathcal{H}^0_Z(\mathcal{C})=0$,
which effectively means that the natural map $\mathcal{C}\longrightarrow j_{\ast
}j^{\ast} \mathcal{C}$ is injective. Hence also the map
$\mathcal{H}om_{\mathcal{O}_Y}(\mathcal{C},\mathcal{C}) \longrightarrow
\mathcal{H}om_{\mathcal{O}_Y} (\mathcal{C}, j_{\ast}j^{\ast} \mathcal{C})\cong
j_{\ast} \mathcal{H}om_{\mathcal{O}_Y}(\mathcal{C} ,\mathcal{C})|_U$ is an injection.
Summing up, we can build a commutative diagram
\[
\begin{CD}
\mathcal{R} @>{\iota_1}>> \mathcal{H}om_{\mathcal{O}_Y} (\mathcal{C},
\mathcal{O}_Y)\cong
\mathcal{H}om_{\mathcal{O}_Y}(\mathcal{C} ,\mathcal{C})\\
@VV{\varepsilon_1}V           @VV{\iota_2}V \\
j_{\ast }\mathcal{R}|_U @>{\varepsilon_2}>>  j_{\ast
}\mathcal{H}om_{\mathcal{O}_Y}(\mathcal{C},
\mathcal{C})|_U 
\end{CD}
\] 
where $\varepsilon_1$ and $\varepsilon_2$ are isomorphisms and $\iota_1$, $\iota_2$
injections ($\iota_1$ is just the natural inclusion of $\mathcal{R}$ into its bidual).
Thus $\iota_1$ and $\iota_2$ are likewise isomorphisms whence $\mathcal{R}\cong
\mathcal{H}om_{\mathcal{O}_Y}(\mathcal{C} ,\mathcal{C})$ and by taking global sections
the desired conclusion $R\cong \mathrm{Hom}_{A_Y}(C,C)$.
\end{proof}
Thus, if one puts the condition ($\heartsuit$), the fact that $R$ is a ring falls
in one's lap almost automatically, a philosophy already supported in [Cat4] (see the
discussion following remark 6.6); but the situation now appears to be even more
idyllic than it was hinted at to be in the latter reference. Note that the crucial
property working behind the scenes and making all the above extension arguments work
is only the Cohen-Macaulayness of $R$.\\

Theorem 1.1 has the following analogue in the realm of graded rings. (In the statement
of the result I will slightly abuse notation falling back on the previously
encountered symbols $R$ etc.; but this is probably more transparent than introducing
too many new letters). 

\begin{theorem}[theorem 1.1 bis]
Let $S=S_0\oplus S_1\oplus\ldots $ be a graded Cohen-Macaulay ring such that $S_0$ is a field and $S$ a finitely generated $S_0-$algebra.  
 Let $R$ be a (graded) finite $S-$module with a minimal graded free resolution
\[
\begin{CD}
0 \longrightarrow \bigoplus_l S(-s_l) @>{\psi }>> \bigoplus_k S(-r_k) @>{\varphi }>>
\bigoplus_j S(-q_j) @>p>> R \longrightarrow 0 .
\end{CD}
\]
Let $\varphi'$ be the submatrix of $\varphi$ consisting of all the rows of $\varphi$
except the first and denote by
$I$ resp.
$I'$ the zeroth Fitting ideals of $\varphi$ resp.
$\varphi'$. \\ Suppose that $\depth (\mathrm{Ann}_S\, R , S)=2$ and
\begin{itemize}
\item[($\heartsuit$)] $\depth \, I' \ge 4$.
\end{itemize}  
Then $R\cong \mathrm{Hom}_{S_Y} (C ,C )$, where $S_Y:=S/\mathrm{Ann}_S\, R$
and $C:=\mathrm{Hom}_{S_Y} (R, S_Y)$, the conductor of $R$ into $S_Y$;
$R$ is given an $S-$algebra structure if we define multiplication as composition
of endomorphisms of $C$.  
\end{theorem}

The proof is completely analogous to that of theorem 1.1. The assumptions on $S$ are made to be able to apply the Auslander-Buchsbaum formula in the graded case to conclude as before that $R$ is Cohen-Macaulay. One then works with the sheaf $\mathcal{R}$ associated to $R$ on $\mathrm{Spec}\, S_Y$, the affine cone over $\mathrm{Proj}\, S_Y$, to finish the argument.

\section{The Gorenstein condition and self-duality of the Hilbert
resolution}
In this section I tie theorem 1.1 in with the symmetry properties that are native to
Gorenstein algebras; in fact, this aspect of the problem has already received
satisfactory treatment in [Gra], but I will outline proofs for the sake of
completeness (also because the arguments in [Gra] are closely interwoven with the
somewhat independent concept of a \emph{Koszul module} and occasionally there is room
for simplifications).\\
The basic set-up will be kept. In particular $A$ is a Cohen-Macaulay local ring and
$R$ a finite $A-$ module. First recall
\begin{definition}
If $R$ is a perfect $A-$ algebra (meaning that $\depth(\mathrm{Ann}_A \, R
,A)=\projdim_A\, R$) and
$R\cong
\mathrm{Ext}^c_A(R,A)$ as
$R-$modules (where 
$c=\dim A-\dim_A R$) then
$R$ is said to be a \emph{Gorenstein} $A-$\emph{algebra of codimension} $c$. 
\end{definition}

\begin{remark}
The $R-$module structure on the a priori $A-$module $\mathrm{Ext}_A^c(R,A)$ is induced
from $R$ by functoriality of $\mathrm{Ext}_A^c(\cdot ,A)$: Thus if for $r\in R$
$\mathrm{mult}_r$ is multiplication by $r$ on $R$, then
$\mathrm{Ext}_A^c(\mathrm{mult}_r ,A)$ is multiplication by $r$ on
$\mathrm{Ext}^c_A(R,A)$.
\end{remark}

Similarly if $S=S_0\oplus S_1\oplus \ldots$ is a positively graded Cohen-Macaulay
ring with $S_0$ a field, and $S$ finitely generated over $S_0$ as an algebra, one can
make
\begin{definition}
A finite perfect graded $S-$algebra $B$ is called a \emph{Gorenstein} $S-$\emph{algebra
of codimension} $c$ (and \emph{with twist} $t\in\mathbb{Z}$) if
$B\cong\mathrm{Ext}^c_S(B, S(t))$ as $B-$modules where $c=\dim S-\dim_S B$.  
\end{definition}
Then we have the following characterization:
\begin{theorem}
If the finite $A-$module $R$ ($A$ CM local and $2$ invertible in
$A$) is a Gorenstein algebra of codimension $2$, then it admits a symmetric minimal
free resolution of the form
\begin{displaymath}
\begin{CD}
\mathbf{R}_{\bullet}:\quad 0\longrightarrow (A^n)^{\vee}
@>{\psi={{-\beta^{\vee}}\choose {\alpha^{\vee}}}}>> A^n\oplus (A^n)^{\vee} @>{
\varphi=(\alpha\, \beta )}>> A^n @>p>> R \longrightarrow 0 .
\end{CD}
\end{displaymath}
Conversely, if $I'$ denotes the zeroth Fitting ideal of the submatrix $\varphi'$ of
$\varphi$ consisting of all the rows of $\varphi$ except the first and if we have 
$\depth I'\ge 4$, then the existence of a symmetric minimal free resolution of the
form $\mathbf{R}_{\bullet }$ is also \emph{sufficient} for the
$A-$module $R$ to be a Gorenstein $A-$algebra of codimension $2$.   
\end{theorem}

\begin{proof}
Let me start with proving the converse: Thus one is given a finite $A-$module $R$
with resolution $\mathbf{R}_{\bullet}$ and $\depth I'\ge 4$; by theorem 1.1 and
remark 1.2 one knows that $R$ carries the structure of a commutative $A-$algebra with
$1$, and moreover $\depth (\mathrm{Ann}_A\, R ,A)=\mathrm{min} \{ i\, :\,
\mathrm{Ext}^i_A(R,A)\neq 0\}=2$ whence I can dualize $\mathbf{R}_{\bullet}$ and
build the following commutative diagram\\
\setlength{\unitlength}{2cm}
\begin{picture}(6,1.7)
\put(0,1){$0$} \put(0.7,1){$(A^n)^{\vee}$} \put(2.2,1){$A^n\oplus (A^n)^{\vee}$}
\put(4,1){$A^n$} \put(5,1){$R$} \put(6,1){$0$} 

\put(0,0){$0$} \put(0.7,0){$(A^n)^{\vee}$} \put(2.2,0){$(A^n)^{\vee}\oplus
A^n$}
\put(4,0){$A^n$} \put(4.7,0){$\mathrm{Ext}^2_A(R,A)$} \put(6,0){$0$} 

\put(0.2,1.05){\vector(1,0){0.4}}  \put(0.2,0.05){\vector(1,0){0.4}}
\put(1.3,1.05){\vector(1,0){0.8}}  \put(1.3,0.05){\vector(1,0){0.8}}
\put(3.25,1.05){\vector(1,0){0.7}} \put(3.25,0.05){\vector(1,0){0.7}}
\put(4.3,1.05){\vector(1,0){0.5}}  \put(4.3,0.05){\vector(1,0){0.3}}
\put(5.25,1.05){\vector(1,0){0.6}} \put(5.7,0.05){\vector(1,0){0.2}}

\put(0.9,0.85){\vector(0,-1){0.6}}  \put(2.6,0.85){\vector(0,-1){0.6}}
\put(4.1,0.85){\vector(0,-1){0.6}}  \put(5.1,0.4){\vector(0,-1){0.15}}
\put(5.1,0.85){\line(0,-1){0.1}} \put(5.1,0.7){\line(0,-1){0.1}}
\put(5.1,0.55){\line(0,-1){0.1}}

\put(3.4,1.1){\small{$(\alpha\,\beta)$}}  \put(1.4,1.15){\small{$-\beta^{\vee} \choose
\alpha^{\vee}$}}
\put(1.5,0.15){\small{$\alpha^{\vee} \choose \beta^{\vee}$}}
\put(3.3,0.1){\small{$(-\beta\, \alpha)$}} 

\put(4.2,0.5){\small{$\mathrm{id}_{A^n}$}}
\put(2.6,0.5){\footnotesize{$\left(\begin{array}{@{}c@{}c@{}}0 &
\mathrm{id}_{(A)^{\vee}}
\\
                                       -\mathrm{id}_{A^n} & 0\end{array}\right)$}}

\put(1,0.5){\small{$\mathrm{id}_{(A^n)^{\vee}}$}}
\put(5.2,0.5){\small{$\cong$}}
\put(4.95,0.5){\small{$u$}}
\end{picture}
\\
\\
\\
where $u: R\cong \mathrm{Ext}^2_A(R,A)$ is the $A-$module isomorphism induced by the
diagram. I claim that this is also an isomorphism of $R-$modules: indeed, from the
proof of theorem 1.1 we know that $R$ is a subring of $A_Y[d^{-1}]$ where
$A_Y:=A/\mathrm{Ann}_A\, R$ and $d\in A_Y$ is some nonzerodivisor. Let then $a/d$,
$a\in A_Y$, be some element of $R$. Then for $r\in R$ we have $u(d(a/d) r)=a u(r)= d
u((a/d)r)$, i.e. $u((a/d)r)= (a/d) u(r)$ and $u$ is also $R-$ linear.

To prove the other direction, let $R$ be a Gorenstein algebra of codimension $2$, and
let an isomorphism $u: R\cong \mathrm{Ext}^2_A(R,A)$ be given. Since $R$ is perfect
and $\dim A-\dim_A R=\depth (\mathrm{Ann}_A\, R ,A)=2$ ($A$ is CM) we can take a
minimal free length
$2$ resolution of
$R$
\[
\begin{CD}
0@>>> F_2 @>{\Psi}>> F_1 @>{\Phi}>> F_0 @>>> R @>>> 0
\end{CD}
\]
where from $R\cong \mathrm{Ext}^2_A(R,A)$ it follows that $\mathrm{rank}\,
F_0=\mathrm{rank}\, F_2$, and since $\mathrm{Ann}_A\, R\neq 0$ $\mathrm{rank}\,
F_0-\mathrm{rank}\, F_1 +\mathrm{rank}\, F_2=0$ whence I can write $F_0\cong A^n$,
$F_1\cong A^n\oplus (A^n)^{\vee}$, $F_2\cong (A^n)^{\vee}$ for some integer $n$.
 Dualize this resolution to obtain a minimal free resolution of
$\mathrm{Ext}^2_A(R,A)$ and lift the given isomorphism $u$ to an isomorphism of
minimal free resolutions:\\
\setlength{\unitlength}{2cm}
\begin{picture}(6,1.6)
\put(0,1){$0$} \put(0.7,1){$(A^n)^{\vee}$} \put(2.2,1){$A^n\oplus (A^n)^{\vee}$}
\put(4,1){$A^n$} \put(5,1){$R$} \put(6,1){$0$} 

\put(0,0){$0$} \put(0.7,0){$(A^n)^{\vee}$} \put(2.2,0){$(A^n)^{\vee}\oplus
A^n$}
\put(4,0){$A^n$} \put(4.7,0){$\mathrm{Ext}^2_A(R,A)$} \put(6,0){$0$} 

\put(0.2,1.05){\vector(1,0){0.4}}  \put(0.2,0.05){\vector(1,0){0.4}}
\put(1.3,1.05){\vector(1,0){0.8}}  \put(1.3,0.05){\vector(1,0){0.8}}
\put(3.25,1.05){\vector(1,0){0.7}} \put(3.25,0.05){\vector(1,0){0.7}}
\put(4.3,1.05){\vector(1,0){0.5}}  \put(4.3,0.05){\vector(1,0){0.3}}
\put(5.25,1.05){\vector(1,0){0.6}} \put(5.7,0.05){\vector(1,0){0.2}}

\put(0.9,0.85){\vector(0,-1){0.6}}  \put(2.6,0.85){\vector(0,-1){0.6}}
\put(4.1,0.85){\vector(0,-1){0.6}}  \put(5.1,0.85){\vector(0,-1){0.6}}

\put(3.5,1.1){\small{$\Phi$}}  \put(1.5,1.1){\small{$\Psi$}}
\put(1.5,0.1){\small{$\Phi^{\vee}$}}
\put(3.5,0.1){\small{$\Psi^{\vee}$}} 

\put(4.2,0.5){\small{$f_1$}}
\put(2.7,0.5){\small{$f_2$}}

\put(1,0.5){\small{$f_3$}}
\put(5.2,0.5){\small{$\cong$}}
\put(4.95,0.5){\small{$u$}}
\end{picture}
\\
\\
\\
(Here and in the following I always implicitly identify free modules and
maps between them with their double duals). Dualizing once more, we can build the
following diagram\\
\setlength{\unitlength}{2cm}
\begin{picture}(6.5,2.6)
\put(0,1){$0$} \put(0.7,1){$(A^n)^{\vee}$} \put(2.2,1){$A^n\oplus (A^n)^{\vee}$}
\put(4,1){$A^n$} \put(4.6,1){\small{$\mathrm{Ext}^2_A(\mathrm{Ext}^2_A(R,A),A)$}}
\put(6.5,1){$0$} 

\put(0,0){$0$} \put(0.7,0){$(A^n)^{\vee}$} \put(2.2,0){$(A^n)^{\vee}\oplus
A^n$}
\put(4,0){$A^n$} \put(4.7,0){$\mathrm{Ext}^2_A(R,A)$} \put(6.5,0){$0$} 

\put(0.2,1.05){\vector(1,0){0.4}}  \put(0.2,0.05){\vector(1,0){0.4}}
\put(1.3,1.05){\vector(1,0){0.8}}  \put(1.3,0.05){\vector(1,0){0.8}}
\put(3.25,1.05){\vector(1,0){0.7}} \put(3.25,0.05){\vector(1,0){0.7}}
\put(4.3,1.05){\vector(1,0){0.2}}  \put(4.3,0.05){\vector(1,0){0.3}}
\put(6.25,1.05){\vector(1,0){0.2}} \put(5.7,0.05){\vector(1,0){0.7}}

\put(0.9,0.85){\vector(0,-1){0.6}}  \put(2.6,0.85){\vector(0,-1){0.6}}
\put(4.1,0.85){\vector(0,-1){0.6}}  \put(5.1,0.85){\vector(0,-1){0.6}}

\put(3.5,1.1){\small{$\Phi$}}  \put(1.5,1.1){\small{$\Psi$}}
\put(1.5,0.1){\small{$\Phi^{\vee}$}}
\put(3.5,0.1){\small{$\Psi^{\vee}$}} 

\put(4.2,0.5){\small{$f_3^{\vee }$}}
\put(2.7,0.5){\small{$f_2^{\vee}$}}

\put(1,0.5){\small{$f_1^{\vee}$}}
\put(5.2,0.5){\small{$\mathrm{Ext}^2_A(u,A)$}}

\put(0,2){$0$} \put(0.7,2){$(A^n)^{\vee}$} \put(2.2,2){$A^n\oplus (A^n)^{\vee}$}
\put(4,2){$A^n$} \put(5,2){$R$} \put(6.5,2){$0$}
 
\put(0.9,1.85){\vector(0,-1){0.6}}  \put(2.6,1.85){\vector(0,-1){0.6}}
\put(4.1,1.85){\vector(0,-1){0.6}}  

\put(0.2,2.05){\vector(1,0){0.4}}  
\put(1.3,2.05){\vector(1,0){0.8}}  
\put(3.25,2.05){\vector(1,0){0.7}} 
\put(4.3,2.05){\vector(1,0){0.5}}  
\put(5.25,2.05){\vector(1,0){1.1}} 

\put(3.5,2.1){\small{$\Phi$}}  \put(1.5,2.1){\small{$\Psi$}}

\put(4.2,1.5){\small{$\mathrm{id}_{A^n}$}}
\put(2.7,1.5){\small{$\mathrm{id}_{A^n\oplus (A^n)^{\vee}}$}}
\put(1,1.5){\small{$\mathrm{id}_{(A^n)^{\vee}}$}}

\put(5.2,1.5){\small{$\cong$}}
\put(4.95,1.5){\small{$c$}}

\put(5.1,1.4){\vector(0,-1){0.15}}
\put(5.1,1.85){\line(0,-1){0.1}} \put(5.1,1.7){\line(0,-1){0.1}}
\put(5.1,1.55){\line(0,-1){0.1}}
\end{picture}
\\
\\
\\
where $c$ is the canonical isomorphism determined by the diagram. Now the point that
needs some work is to compare the isomorphisms $u$ and $\mathrm{Ext}^2_A(u,A)\circ c$.
This is done in [Gra] (cf. especially lemma 2.1 and thm. 3) by identifying the
functors $\mathrm{Ext}_A^2(\cdot ,A)$ resp. $\mathrm{Ext}^2_A(\mathrm{Ext}^2_A(\cdot
,A) , A)$ with $\mathrm{Hom}_A(\cdot , A/(x_1, x_2))$ resp.
$\mathrm{Hom}_A(\mathrm{Hom}_A(\cdot ,A/(x_1 ,x_2)) , A/(x_1 ,x_2))$ where $x_1 ,\:
x_2$ is a regular sequence in $\mathrm{Ann}_A\, R$, and making all of the occurring
isomorphisms explicit. Grassi uses the extra assumption that $A$ be a domain to make the proof work, but this is in fact redundant as I will show in lemma 3.2 below (the argument needed to remove this hypothesis is a little technical, so I deferred it to the next section where, after all, it integrates rather better). It
turns out that $\mathrm{Ext}^2_A(u,A)\circ c=-u$, and since $2$ was supposed to be
invertible in $A$, one can conclude that the following diagram commutes:\\
\setlength{\unitlength}{2cm}
\begin{picture}(6,1.6)
\put(0,1){$0$} \put(0.7,1){$(A^n)^{\vee}$} \put(2.2,1){$A^n\oplus (A^n)^{\vee}$}
\put(4,1){$A^n$} \put(5,1){$R$} \put(6,1){$0$} 

\put(0,0){$0$} \put(0.7,0){$(A^n)^{\vee}$} \put(2.2,0){$(A^n)^{\vee}\oplus
A^n$}
\put(4,0){$A^n$} \put(4.7,0){$\mathrm{Ext}^2_A(R,A)$} \put(6,0){$0$} 

\put(0.2,1.05){\vector(1,0){0.4}}  \put(0.2,0.05){\vector(1,0){0.4}}
\put(1.3,1.05){\vector(1,0){0.8}}  \put(1.3,0.05){\vector(1,0){0.8}}
\put(3.25,1.05){\vector(1,0){0.7}} \put(3.25,0.05){\vector(1,0){0.7}}
\put(4.3,1.05){\vector(1,0){0.5}}  \put(4.3,0.05){\vector(1,0){0.3}}
\put(5.25,1.05){\vector(1,0){0.6}} \put(5.7,0.05){\vector(1,0){0.2}}

\put(0.9,0.85){\vector(0,-1){0.6}}  \put(2.6,0.85){\vector(0,-1){0.6}}
\put(4.1,0.85){\vector(0,-1){0.6}}  \put(5.1,0.85){\vector(0,-1){0.6}}

\put(3.5,1.1){\small{$\Phi$}}  \put(1.5,1.1){\small{$\Psi$}}
\put(1.5,0.1){\small{$\Phi^{\vee}$}}
\put(3.5,0.1){\small{$\Psi^{\vee}$}} 

\put(4.2,0.5){\small{$\frac{f_1-f_3^{\vee}}{2}$}}
\put(2.7,0.5){\small{$\frac{f_2-f_2^{\vee}}{2}$}}

\put(1,0.5){\small{$\frac{f_3-f_1^{\vee}}{2}$}}
\put(5.2,0.5){\small{$\cong$}}
\put(4.95,0.5){\small{$u$}}
\end{picture}
\\
\\
\\
Then $(f_2-f_2^{\vee})/2$ is a skew isomorphism which, by a suitable orthogonal
isomorphism $B$ of 
$A^n\oplus (A^n)^{\vee}$, $B^{\vee}B=\mathrm{id}_{A^n\oplus (A^n)^{\vee}}$, can be
brought to normal form
\[
J:=\left(\begin{array}{@{}c@{}c@{}}0 &
\mathrm{id}_{(A)^{\vee}}
\\
                                       -\mathrm{id}_{A^n} & 0\end{array}\right)
=B^{\vee} \left( \frac{f_2-f_2^{\vee}}{2} \right) B .
\]
One is done because the lower row of the following commutative diagram is of the form
required in resolution $\mathbf{R}_{\bullet}$:\\
\setlength{\unitlength}{2cm}
\begin{picture}(6,1.6)
\put(0,1){$0$} \put(0.7,1){$(A^n)^{\vee}$} \put(2.2,1){$A^n\oplus (A^n)^{\vee}$}
\put(4,1){$A^n$} \put(5,1){$R$} \put(6,1){$0$} 

\put(0,0){$0$} \put(0.7,0){$(A^n)^{\vee}$} \put(2.2,0){$A^n\oplus
(A^n)^{\vee}$}
\put(4,0){$A^n$} \put(4.7,0){$\mathrm{Ext}^2_A(R,A)$} \put(6,0){$0$} 

\put(0.2,1.05){\vector(1,0){0.4}}  \put(0.2,0.05){\vector(1,0){0.4}}
\put(1.3,1.05){\vector(1,0){0.8}}  \put(1.3,0.05){\vector(1,0){0.8}}
\put(3.25,1.05){\vector(1,0){0.7}} \put(3.25,0.05){\vector(1,0){0.7}}
\put(4.3,1.05){\vector(1,0){0.5}}  \put(4.3,0.05){\vector(1,0){0.3}}
\put(5.25,1.05){\vector(1,0){0.6}} \put(5.7,0.05){\vector(1,0){0.2}}

\put(0.9,0.85){\vector(0,-1){0.6}}  \put(2.6,0.85){\vector(0,-1){0.6}}
\put(4.1,0.85){\vector(0,-1){0.6}}  \put(5.1,0.85){\vector(0,-1){0.6}}

\put(3.5,1.1){\small{$\Phi$}}  \put(1.5,1.1){\small{$\Psi$}}
\put(1.25,0.1){\small{$-J (\Phi\circ B)^{\vee}$}}
\put(3.4,0.1){\small{$\Phi\circ B$}} 

\put(4.2,0.5){\small{$\mathrm{id}_{A^n}$}}
\put(2.7,0.5){\small{$B^{\vee}$}}

\put(1,0.5){\small{$\frac{f_3-f_1^{\vee}}{2}$}}
\put(5.2,0.5){\small{$\mathrm{id}_R$}}

\end{picture}
\\
\\
\\
\end{proof}

Finally let me state the theorem corresponding to theorem 2.4 in the graded case
whose proof is entirely similar to the preceding one.

\begin{theorem}[theorem 2.4 bis]
Let $S=S_0\oplus S_1\oplus\ldots $ be a graded Cohen-Macaulay ring 
such that $S_0$ is a field and $S$ a finitely generated $S_0-$algebra. Assume
$\mathrm{char}\, S_0\neq 2$. Then every Gorenstein $S-$algebra of codimension $2$ and
with twist $t$ has a symmetric graded free resolution 
\begin{eqnarray*}
\begin{CD}
\mathbf{R}_{\bullet}:\; 0\longrightarrow \bigoplus_{j=1}^nS^{\vee}(-t+r_j)
@>{\psi={{-\beta^{\vee}}\choose {\alpha^{\vee}}}}>> \bigoplus_{k=1}^n
S^{\vee}(-t+s_k)\oplus \bigoplus_{k=1}^n S(-s_k)
\end{CD}\\
\begin{CD}
 @>{
\varphi=(\alpha\, \beta )}>> \bigoplus_{j=1}^n S(-r_j) @>p>> R \longrightarrow 0 ,
\end{CD}\quad
\end{eqnarray*}
($n\in\mathbb{N},\; t,\: r_j ,\: s_k\in\mathbb{Z}$). Conversely, if $I'$ denotes the
zeroth Fitting ideal of the matrix $\varphi'$ which is the matrix $\varphi$ with first
row erased and 
$\depth I'\ge 4$, then a graded $S-$module $R$ with a symmetric minimal
graded free resolution of the form $\mathbf{R}_{\bullet }$ is a 
Gorenstein $S-$algebra of codimension $2$ and with twist $t$.  \flushright{\qedsymbol} 
\end{theorem}
The proof is almost entirely similar to the preceding one, but there is one point that deserves mentioning: When proving that the Gorenstein $S-$algebra $R$ has a symmetric graded free resolution, I want to use that the given isomorphism $R\cong \mathrm{Ext}^2_S (R, S(t))$ is skew-symmetric with respect to the duality given by $\mathrm{Ext}^2_S(- , S(t))$. In the local case this followed from [Gra]. But then, in particular, in the present situation, one knows that for every $\mathfrak{p}\in\mathrm{Supp}\, R$ the localized isomorphism $R_{\mathfrak{p}}\cong \mathrm{Ext}_{S_{\mathfrak{p}}}(R_{\mathfrak{p}}, S_{\mathfrak{p}})$ is skew-symmetric, whence also the original isomorphism $R\cong \mathrm{Ext}^2_S(R, S(t))$ is skew-symmetric. 

\section{Regularizing minors}
\setcounter{equation}{0}
This section provides a strengthening of theorem 2.4, the local case of the strucuture theorem, and thereby clarifies the relationship between Gorenstein algebras in codimension 2 and Koszul modules as introduced in [Gra].\\
Therefore let $(A, \mathfrak{m} ,k)$ be again CM local with $2\notin \mathfrak{m}$ and $R$ a codimension 2 Gorenstein algebra over $A$. 
Whereas the usual Koszul complex is associated with a linear form $f: A^n\to A$, a
Koszul module is a module having a resolution similar to the Koszul complex up to
the fact that the r\^{o}le of $f$ is taken by a family of (vector-valued) maps from
$A^n$ to $A^n$. I'll only make this precise in the relevant special case:
\begin{definition}
A finite $A-$module $M$ having a length 2 resolution 
\begin{equation}
\begin{CD}
0 \to A^n @>{\rho_1\choose\rho_2}>> A^{2n} @>(\tau_1\:\tau_2)>> A^n \to M \to 0
\end{CD}
\end{equation}
some $n\in\mathbb{N}$, is a Koszul module iff
$\mathrm{det}(\tau_1),\;\mathrm{det}(\tau_2)$ is a regular sequence on $A$ and
$\exists$ a unit
$\lambda\in A$:
$\mathrm{det}(\rho_1)=(-1)^n\lambda
\mathrm{det}(\tau_2),\;\mathrm{det}(\rho_2)=\lambda \mathrm{det}(\tau_1)$.
\end{definition}
Then Grassi proves in case $A$ is a domain ([Gra], thm. 3.3) that $R$ has a
(Gorenstein) symmetric resolution
\begin{equation}
\begin{CD}
0 \to A^n @>{-\beta^t\choose\alpha^t}>> A^{2n} @>(\alpha\:\beta)>> A^n \to R \to 0
\end{CD}
\end{equation}
and a second resolution of the prescribed type (1) for the Koszul module
condition, and that these 2 are related by an isomorphism of complexes which is
the identity in degrees 0 and 2; firstly, for sake of generality, I will briefly show
that the assumption "$A$ a domain" is in fact not needed, thereby closing also the remaining gap in the proof of theorem 2.4, and secondly, prove that
there is one single resolution of $R$ meeting both requirements, i.e. a resolution
as in (2) with $\mathrm{det}(\alpha),\;\mathrm{det}(\beta)$ an $A-$regular
sequence.
\begin{lemma}
A Gorenstein algebra $R$ has a minimal free resolution of type (2) over any CM local \emph{ring} with $2\notin\mathfrak{m}$ (i.e. one need not assume that $A$ be a domain).
\end{lemma}
\begin{proof}
Note that the only place in [Gra] where the hypothesis that
$A$ be a domain enters is at the beginning of the proof of proposition 1.5, page
930: Here one is given a resolution as in (1), but \underline{without} any
additional assumptions on
$\mathrm{det}(\tau_1),\mathrm{det}(\tau_2),\mathrm{det}(\rho_1),\mathrm{det}(\rho_2)$
whatsoever, and Grassi wants to conclude that $\exists$ a base change in $A^{2n}$
such that (in the new base) $\mathrm{det}(\tau_1)$ is not a zero divisor on $A$.
But this can be proven by a similar method as Grassi uses in the sequel of the
proof of proposition 1.5, without using "$A$ a domain":
For let $\mathfrak{p}_1,\ldots,\mathfrak{p}_r$ be the associated primes of $A$
which are precisely the minimal elements of $\mathrm{Spec}(A)$ since $A$ is CM.
One shows that $\exists$ a base change in $A^{2n}$ such that
$\mathrm{det}(\tau_1)\notin\mathfrak{p}_i,\;\forall i=1,\ldots,r$ (in the new
base), more precisely, that $\exists$ a sequence of $r$ base changes such that
after the $m$th base change
\[
(\ast) \hspace{1cm} \mathrm{det}(\tau_1)\notin\mathfrak{p}_i,\;\forall
i\in\{r-m+1,\ldots,r\},
\] 
$m=0,\ldots,r$, the assertion being empty for $m=0$. Therefore, inductively,
suppose $(\ast)$ holds for $m$ to get it for $m+1$.\\
Denote by $[i_1,\ldots,i_n]$ the maximal minor of $(\tau_1\;\tau_2)$
corresponding to the columns $i_1,\ldots,i_n,\;i_j\in\{1,\ldots,2n\}$. If
$[1,\ldots,n]\notin\mathfrak{p}_{r-m}$ I'm already O.K., so suppose
$[1,\ldots,n]\in\mathfrak{p}_{r-m}$. By the Eisenbud-Buchsbaum acyclicity
criterion $I_n((\tau_1\;\tau_2))$ cannot consist of zerodivisors on $A$ alone,
therefore set 
\begin{eqnarray*}
l_1:=\mathrm{min}\{c :\:\exists s_1,\ldots,s_{n-1}\:\mathrm{with}\:
s_1<s_2<\ldots <s_{n-1}<c \\ \mathrm{and}  
 [s_1,\ldots ,s_{n-1},c]\notin\mathfrak{p}_{r-m}\} \end{eqnarray*} (then $n<l_1\le
2n$) and inductively,
\begin{eqnarray*}
l_i:=\mathrm{min}\{c :\:\exists s_1',\ldots,s_{n-i}'\:\mathrm{with}\:
s_1'<\ldots<s_{n-i}'<c<l_{i-1}<\ldots<l_1\\
\mathrm{and}\:[s_1',\ldots,s_{n-i}',c,l_{i-1},\ldots,l_1]\notin\mathfrak{p}_{r-m}\}
,\end{eqnarray*} 
$i=2,\ldots,n$. Then $\exists \;J$ such that $n<l_J<l_{J-1}<\ldots<l_1\le 2n$ and
for $I>J\;l_I\in\{1,\ldots,n\}$ ($J=n$ might occur and then the set of
$l_I\in\{1,\ldots,n\}$ is empty; this does not matter). \\
I have $[l_n,\ldots,l_1]\notin\mathfrak{p}_{r-m}$ by construction. Choose
$b\in(\bigcap\limits_{i=r-m+1}^{r}\mathfrak{p}_i)\backslash\mathfrak{p}_{r-m}$,
which is nonempty since the $\mathfrak{p}_i$'s are the minimal elements of
$\mathrm{Spec}(A)$. Denote by $y_1< \ldots <y_J$ the complementary indices of the
$l_I\in\{1,\ldots ,n\}$ inside $\{ 1,\ldots ,n \}$ and consider the base change on
$A^{2n}$: $M_{y_1,l_J}(b)\circ M_{y_2,l_{j-1}}(b)\circ \ldots \circ
M_{y_J,l_1}(b)$, where $M_{y_{\nu},l_{J-\nu +1}}(b),\; \nu=1,\ldots ,J$ is addition
of $b$ times the $l_{J-\nu +1}$ column to the $y_{\nu}$ column. Then one sees (by
the multilinearity of determinants) 
\[
[1,\ldots,n]_{\mathrm{new}}=[1,\ldots,n]_{\mathrm{old}}\pm
b^J[l_n,\ldots,l_1]_{\mathrm{old}}+b\mu,
\]
where "new" means after and "old" before the base change and $\mu$ is an
element in $\mathfrak{p}_{r-m}$ by the defining minimality property of the $l$'s.
Therefore, since by the induction hypothesis $[1,\ldots,n]_{\mathrm{old}}\notin
\mathfrak{p}_i,\; \forall i\in\{r-m+1,\ldots,r\}$ and $b$ is chosen appropriately:
$[1,\ldots,n]_{\mathrm{new}}\notin \mathfrak{p}_i,\; \forall
i\in\{r-m,\ldots,r\}$. This finally proves $\mathrm{det}(\tau_1) \notin
\mathfrak{p}_i \; \forall i=1,\ldots,r$ after the sequence of base changes, i.e.
$\mathrm{det}(\tau_1)$ is then $A$-regular, that what was to be shown.\end{proof}
\vspace{0.5cm}
Secondly, I now want to prove:
\begin{theorem}
A codimension 2 Gorenstein algebra $R$ over a local CM ring $(A,\mathfrak{m},k)$ with
$2\notin \mathfrak{m}$ has a resolution
\[
\begin{CD}
0 \to A^n @>{-\beta^t\choose\alpha^t}>> A^{2n} @>(\alpha\:\beta)>> A^n \to R \to 0
\end{CD}
\]
which is also of Koszul module type, i.e.
$\mathrm{det}(\alpha),\:\mathrm{det}(\beta)$ is an $A-$regular sequence.
\end{theorem}
\begin{proof}
Taking into account the above remark that one can dispose of the assumption "$A$ a
domain" the fact that $R$ has a resolution with the symmetry property above is proven
in [Gra], thm. 3.3., so I have to show that $\exists$ a base change in $A^{2n}$ which
preserves the relation $\alpha \beta^t=\beta \alpha^t$ and in the new base
$\mathrm{det}(\alpha),\:\mathrm{det}(\beta)$ is a regular sequence. The punch line
to show this is as in the foregoing argument except that everything is a little
harder because one has to keep track of preserving the symmetry: Therefore let again
be  $\mathfrak{p}_1,\ldots,\mathfrak{p}_r$ the associated primes of $A$, and I show
that $\exists$ a sequence of $r$ base changes in $A^{2n}$ preserving the symmetry
and such that after the $m$th base change $(\ast)$ above holds, the case $m=0$ being
trivial. For the inductive step, suppose $\det(\alpha)\in\mathfrak{p}_{r-m}$ to
rule out a trivial case; I write
$[i_1,\ldots,i_{\nu};j_1,\ldots,j_{n-
\nu}]\equiv\det(\alpha_{i_1}\ldots\alpha_{i_{\nu}}\:\beta_{
j_1}\ldots\beta_{j_{n- \nu}})$. Call a minor $[i_1,\ldots,i_{\nu};j_1,\ldots,j_{n-
\nu}]$
\underline{good} iff $\{ i_1,\ldots,i_{\nu}\} \cap \{j_1,\ldots,j_{n- \nu}
\}=\emptyset$.\\ I want to find a good minor that does not belong to
$\mathfrak{p}_{r-m}$ (possibly after a base change in $A^{2n}$). Therefore suppose
all the good minors belong to
$\mathfrak{p}_{r-m}$. Since $\mathrm{grade}\:I_n((\alpha\:\beta))\ge 2$ by
Eisenbud-Buchsbaum acyclicity, $\exists$ a minor $\notin \mathfrak{p}_{r-m}$
(which is not good). For $n=1$ this is a contradiction since all minors are good, and
I can suppose $n>1$ in the process of finding a good minor. Now choose a minor
$[I_1,\ldots,I_k;J_1,\ldots,J_{n-k}]$ such that
\begin{itemize}
\item
$[I_1,\ldots,I_k;J_1,\ldots,J_{n-k}]\notin \mathfrak{p}_{r-m}$
\item
$\mathrm{card}( \{ I_1,\ldots,I_k \} \cap \{J_1,\ldots,J_{n-k} \} )=:M_0$  is minimal
among the minors which do not belong to $\mathfrak{p}_{r-m}$.
\end{itemize}
I want to perform a base change in $A^{2n}$ not destroying the symmetry such that in
the new base
$\exists$ a minor $[T_1,\ldots,T_{k-1};S_1,\ldots,S_{n-k+1}]$ such that
\begin{itemize}
\item
$[T_1,\ldots,T_{k-1};S_1,\ldots,S_{n-k+1}]\notin \mathfrak{p}_{r-m}$
\item
$\mathrm{card}( \{ T_1,\ldots,T_{k-1} \} \cap \{ S_1,\ldots,S_{n-k+1} \} )=M_0-1.$
\end{itemize}
Continuing this process $M_0$ steps (i.e. performing $M_0$ successive base changes) I
can find a good minor not contained in $\mathfrak{p}_{r-m}$.\\
Let now $[T_1,\ldots,T_{k-1};S_1,\ldots,S_{n-k+1}]$ be given. Choose $H \in \{
I_1,\ldots,I_k \} \cap \{ J_1,\ldots,J_{n-k} \}$ and $L \in \{ 1,\ldots,n \}
- \{ I_1,\ldots,I_k \} \cup \{ J_1,\ldots,J_{n-k} \}$ (both of which exist).
Now perform the base change in $A^{2n}$ which corresponds to adding $\alpha_H$ to
$\beta_L$ and $\alpha_L$ to $\beta_H$ (preserving the symmetry), and consider
\[
\det(\alpha_{I_1}\ldots\hat{\alpha_H}\ldots\alpha_{I_k}\;\beta_{J_1}\ldots\beta_H+\alpha_L
\ldots\beta_L+\alpha_H\ldots\beta_{J_{n-k}}),
\]
an $n\times n-$minor of the transformed matrix which I can write as\\
$[T_1,\ldots,T_{k-1};S_1,\ldots,S_{n-k+1}]$, where $\{ T_1,\ldots,T_{k-1} \} =\{
I_1,\ldots,I_k \}-\{ H\}$, $\{ S_1,\ldots,S_{n-k+1} \} = \{J_1,\ldots,J_{n-k} \}
\cup \{ L \}$ and obviously, $\mathrm{card}( \{ T_1,\ldots,T_{k-1} \}$\\ $\cap \{
S_1,\ldots,S_{n-k+1} \} )=M_0-1.$ I want to prove that this minor does not belong
to $\mathfrak{p}_{r-m}$. For this I show that in fact
\begin{eqnarray*}
[T_1,\ldots,T_{k-1};S_1,\ldots,S_{n-k+1}]=\pm
[I_1,\ldots,I_k;J_1,\ldots,J_{n-k}] \\+\mathrm{"residual\; terms"},
\end{eqnarray*}
where "residual terms"$\in \mathfrak{p}_{r-m}$. Using the additivity of the
determinant in each column I find that "residual terms" consists of 3 summands two of
which clearly belong to $\mathfrak{p}_{r-m}$ because
$[I_1,\ldots,I_k;J_1,\ldots,J_{n-k}]$ was chosen such that $\mathrm{card}( \{ I_1,\ldots,I_k \} \cap \{J_1,\ldots,J_{n-k} \}
)=:M_0$  was minimal among the minors of the matrix before the base change which did
not belong to
$\mathfrak{p}_{r-m}$, whereas the third summand is (up to sign)
\[
\det(\alpha_{I_1}\ldots\hat{\alpha_H}\ldots\alpha_{I_k}\:\alpha_L\:\beta_{J_1}\ldots
\hat{\beta_H}
\ldots\beta_{J_{n-k}}\:\beta_L).
\]
To show that the latter is in $\mathfrak{p}_{r-m}$ I apply the so-called
"Pl\"ucker relations":\vspace{0.3cm}

\itshape
Given an $M\times N-$matrix, $M\le 
 N,\; a_1,\ldots ,a_p,b_q,\ldots ,b_M,c_1,\ldots ,c_s\in \{ 1,\ldots ,N \} 
 ,\; s=M-p+q-1>M,\; t=M-p > 0 $, one has
\[
(P) \sum\limits_{ {i_1<\ldots <i_t \atop i_{t+1}<\ldots <i_s} \atop \{
i_1,\ldots ,i_s\} =\{ 1,\ldots ,s\} }
\sigma(i_1,\ldots,i_s)[a_1,\ldots ,a_p\: c_{i_1}\ldots c_{i_t}][c_{i_{t+1}}\ldots
c_{i_s}\: b_q\ldots b_M]=0
\]
where $\sigma(i_1,\ldots,i_s)$ is the sign of the permutation ${1,\ldots ,s \choose
i_1,\ldots ,i_s}$ (see e.g. [B-He], lemma 7.2.3, p. 308).
\vspace{0.3cm}

\upshape
In my situation, I let $M:=n,\;N:=2n,\;p:=n-2,\;q:=n+1,\;s:=n+1$ and for the
columns corresponding to the $a$'s above I choose the $n-2$ columns
\[
\alpha_{I_1},\: \alpha_{I_2},\ldots ,\hat{\alpha_H},\ldots ,\alpha_{I_k},\: 
\beta_{J_1}, \ldots ,\hat{\beta_H},\ldots ,\beta_{J_{n-k}} 
\]
(in this order), for the columns corresponding to the $b$'s I choose the empty set
(which is allowable here), and finally for the columns corresponding to the $c$'s the
$n+2$ columns
\begin{displaymath} 
\alpha_H,\: \beta_H,\: \alpha_L,\: \beta_L,\: \alpha_{I_1},\: \alpha_{I_2},\ldots
,\: \hat{\alpha_H},\ldots,\: \alpha_{I_k},\: \beta_{J_1},\ldots ,\:
\hat{\beta_H},\ldots ,\beta_{J_{n-k}} 
\end{displaymath}\\
Applying (P) one gets 6 nonvanishing summands, 4 of which (namely 
\begin{eqnarray*}\det(\alpha_{I_1}\ldots\hat{\alpha_H}\ldots\alpha_{I_k}\:\beta_{J_1}\ldots
\hat{\beta_H}
\ldots\beta_{J_{n-k}}\: \alpha_H\:\alpha_L)\cdot \mathrm{(a \; second\; factor)},\\
\det(\alpha_{I_1}\ldots\hat{\alpha_H}\ldots\alpha_{I_k}\:\beta_{J_1}\ldots
\hat{\beta_H}
\ldots\beta_{J_{n-k}}\: \alpha_H\:\beta_L)\cdot \mathrm{(a \; second\; factor)},\\
\det(\alpha_{I_1}\ldots\hat{\alpha_H}\ldots\alpha_{I_k}\:\beta_{J_1}\ldots
\hat{\beta_H}
\ldots\beta_{J_{n-k}}\: \beta_H\:\alpha_L)\cdot \mathrm{(a \; second\; factor)},\\
\det(\alpha_{I_1}\ldots\hat{\alpha_H}\ldots\alpha_{I_k}\:\beta_{J_1}\ldots
\hat{\beta_H}
\ldots\beta_{J_{n-k}}\: \beta_H\:\beta_L)\cdot \mathrm{(a \; second\;
factor)} ) \end{eqnarray*} 
 are in
$\mathfrak{p}_{r-m}$ by the defining minimality property of $[I_1,\ldots
,I_k;J_1,\ldots  ,J_{n-k}]$ above. The remaining 2 summands add up to (watch the
signs!)\\
$\pm
2 \det(\alpha_{I_1}\ldots\hat{\alpha_H}\ldots\alpha_{I_k}\:\alpha_L\:\beta_{J_1}\ldots
\hat{\beta_H}
\ldots\beta_{J_{n-k}}\:\beta_L) \cdot [I_1,\ldots ,I_k;J_1,\ldots  ,J_{n-k}]$
which therefore is also in $\mathfrak{p}_{r-m}$. But $[I_1,\ldots ,I_k;J_1,\ldots 
,J_{n-k}]\notin \mathfrak{p}_{r-m}$ and 2 is a unit in $A$, therefore
$\det(\alpha_{I_1}\ldots\hat{\alpha_H}\ldots\alpha_{I_k}\:\alpha_L\:\beta_{J_1}\ldots
\hat{\beta_H}
\ldots\beta_{J_{n-k}}\:\beta_L)\in\mathfrak{p}_{r-m}$ as desired, since
$\mathfrak{p}_{r-m}$ is prime.\\
Hence inductively, after $M_0$ base changes in $A^{2n}$, I can find a good minor of
the transformed matrix that is not in $\mathfrak{p}_{r-m}$. I assume $[1,\ldots,n]\in
\mathfrak{p}_{r-m}$. I can now define
\begin{eqnarray*}
l_1:=\mathrm{min}\{c :\:\exists s_1,\ldots,s_{n-1}\:\mathrm{with}\:
s_1<s_2<\ldots<s_{n-1}<c\\
\mathrm{and}\:[s_1,\ldots,s_{n-1},c]\notin\mathfrak{p}_{r-m} 
\mathrm{\underline{and}}\: [s_1,\ldots,s_{n-1},c]\: \mathrm{is \:
good} \}
\end{eqnarray*}
(then $n<l_1\le 2n$) and inductively,
\begin{eqnarray*}
l_i:=\mathrm{min}\{c :\:\exists s_1',\ldots,s_{n-i}'\:\mathrm{with}\:
s_1'<\ldots<s_{n-i}'<c<l_{i-1}<\ldots<l_1\\ 
\mathrm{and}\: [s_1',\ldots,s_{n-i}',c,l_{i-1},\ldots,l_1] \: \mathrm{is
\: good} \\
\mathrm{and}\:[s_1',\ldots,s_{n-i}',c,l_{i-1},\ldots,l_1]\notin\mathfrak{p}_{r-m}\}
.\end{eqnarray*} 
Then $[l_n,\ldots,l_1]\notin \mathfrak{p}_{r-m}$ which is good and can therefore be
written as
$[l_n,\ldots,l_1]=[l^{\alpha}_1,\ldots,l^{\alpha}_h;l^{\beta}_1,\ldots,l^{\beta}_{n-h}]$ 
with $\{ l^{\alpha}_1,\ldots,l^{\alpha}_h\} \cap \{
l^{\beta}_1,\ldots,l^{\beta}_{n-h} \} =\emptyset$. Choose
$b\in(\bigcap\limits_{i=r-m+1}^{r}\mathfrak{p}_i)\backslash\mathfrak{p}_{r-m}$ and
perform a base change in $A^{2n}$ (preserving the symmetry) by adding $b$ times the
$l_i^{\beta}$ column of $\beta$ to the $l_i^{\beta}$ column of $\alpha$, for
$i=1,\ldots,n-h$. Then
\[
[1,\ldots,n]_{\mathrm{new}}=[1,\ldots,n]_{\mathrm{old}}\pm
b^{n-h}[l_n,\ldots,l_1]_{\mathrm{old}}+b\mu,
\]
where $\mu\in\mathfrak{p}_{r-m}$ by the defining minimality property of the
$l$'s. Thus $[1,\ldots,n]_{\mathrm{new}}\notin\mathfrak{p}_i$ for
$i=r-m,\ldots,r$, which is the inductive step for the property ($\ast$).
Therefore after a sequence of base changes that preserve the symmetry
$\alpha\beta^t=\beta\alpha^t$, $\det(\alpha)$ can be made an $A-$regular
element.

Let's sum up: I have that $\det(\alpha)$ is a nonzerodivisor in $A$, and want
to prove that $\exists$ a base change in $A^{2n}$ preserving the symmetry and
leaving $\alpha$ unchanged (i.e. fixing the first $n$ basis vectors of
$A^{2n}$) such that in the new base $\det(\beta)$ is a nonzerodivisor in
$A/(\det(\alpha))$. The argument is almost identical to the preceding one. In
fact, let $\mathfrak{q}_1,\ldots ,\mathfrak{q}_s$ be the associated primes of
$A/(\det(\alpha))$ which are exactly the minimal prime ideals containing
$(\det(\alpha))$ because $A/(\det(\alpha))$ is CM ($A$ is CM and
$\det(\alpha)$ is $A-$regular). Then the part of the above proof starting with
"$\ldots$ the symmetry: Therefore let again
be  $\mathfrak{p}_1,\ldots,\mathfrak{p}_r$ the associated primes of $A$, and I show
that $\exists$ a sequence of $r$ base changes in $A^{2n}$ $\ldots$" and ending
with "$\ldots$ Choose $H \in \{
I_1,\ldots,I_k \} \cap \{ J_1,\ldots,J_{n-k} \}$ and $L \in \{ 1,\ldots,n \}
- \{ I_1,\ldots,I_k \} \cup \{ J_1,\ldots,J_{n-k} \}$ $\ldots$" goes through
\underline{verbatim} (and has to be inserted here) if throughout one replaces
$r$ with $s$, $\det(\alpha)$ with $\det(\beta)$, and the symbol
"$\mathfrak{p}$" with "$\mathfrak{q}$". Thereafter, a slight change is
necessary because in the process of finding a good minor, i.e. in the course
of the $M_0$ base changes on $A^{2n}$ that transform $(\alpha\:\beta)$ s.t. in
the new base $\exists$ a good minor, the shape of $\beta$ is changed. This
change must preserve the property $\det(\beta)\notin\mathfrak{q}_1,\ldots
,\mathfrak{q}_{s-m+1}$ in order not to destroy the induction hypothesis. The
way out is as follows:\\
Choose $\zeta
\in (\bigcap\limits_{i=r-m+1}^{r}\mathfrak{q}_i)\backslash\mathfrak{q}_{r-m}$,
which is possible since the $q$'s all have height 1. Now perform the base change in
$A^{2n}$ which corresponds to adding
$\zeta\alpha_H$ to
$\beta_L$ and $\zeta\alpha_L$ to $\beta_H$ (preserving the symmetry), and consider
\[
\det(\alpha_{I_1}\ldots\hat{\alpha_H}\ldots\alpha_{I_k}\;\beta_{J_1}\ldots\beta_H+\zeta\alpha_L
\ldots\beta_L+\zeta\alpha_H\ldots\beta_{J_{n-k}}),
\]
an $n\times n-$minor of the transformed matrix which I can write as\\
$[T_1,\ldots,T_{k-1};S_1,\ldots,S_{n-k+1}]$, where $\{ T_1,\ldots,T_{k-1} \} =\{
I_1,\ldots,I_k \}-\{ H\}$, $\{ S_1,\ldots,S_{n-k+1} \} = \{J_1,\ldots,J_{n-k} \}
\cup \{ L \}$ and obviously, $\mathrm{card}( \{ T_1,\ldots,T_{k-1} \}$\\ $\cap \{
S_1,\ldots,S_{n-k+1} \} )=M_0-1.$ I want to prove that this minor does not belong
to $\mathfrak{q}_{s-m}$ and furthermore that 
\[
\det(\beta_1\ldots \beta_H+\zeta\alpha_L \ldots \beta_L+\zeta\alpha_H \ldots \beta_n) \notin\mathfrak{q}_1,\ldots
,\mathfrak{q}_{s-m+1}.
\] 
The latter statement is obvious by the choice of $\zeta$ (and multilinearity of
determinants). The former one follows if I show
\begin{eqnarray*}
[T_1,\ldots,T_{k-1};S_1,\ldots,S_{n-k+1}]=\pm
\zeta [I_1,\ldots,I_k;J_1,\ldots,J_{n-k}]\\ +\mathrm{"residual\; terms"},
\end{eqnarray*}
where "residual terms"$\in \mathfrak{q}_{s-m}$ because $\zeta$ and
$[I_1,\ldots ,I_k;J_1,\ldots ,J_{n-k}]$ are both $\notin\mathfrak{q}_{s-m}$ by
assumption. Again "residual terms" consists of 3 summands two of which belong
to $\mathfrak{q}_{s-m}$ because of the defining minimality property of
$[I_1,\ldots,I_k;J_1,\ldots,J_{n-k}]$. The third summand is up to sign 
\[
\zeta \det(\alpha_{I_1}\ldots\hat{\alpha_H}\ldots\alpha_{I_k}\:\alpha_L\:\beta_{J_1}\ldots
\hat{\beta_H}
\ldots\beta_{J_{n-k}}\:\beta_L),
\] 
therefore it suffices to show
$\det(\alpha_{I_1}\ldots\hat{\alpha_H}\ldots\alpha_{I_k}\:\alpha_L\:\beta_{J_1}\ldots
\hat{\beta_H}
\ldots\beta_{J_{n-k}}\:\beta_L)\in \mathfrak{q}_{s-m}$. This is done word by word
as in the passage of the first part of this proof starting with "$\ldots$ I apply
the so-called "Pl\"ucker relations":$\ldots$ " and ending "$\ldots$ since
$\mathfrak{p}_{r-m}$ is prime.$\ldots$ ", taking into account the afore-mentioned
changes in notation.\\
The rest of the proof is as follows: Inductively, I can find a good minor of the
transformed matrix that is not in $\mathfrak{q}_{s-m}$. To avoid a trivial case,
I assume $[n+1,\ldots ,2n]\in \mathfrak{q}_{s-m}$. Now I define
\begin{eqnarray*}
L_1:=\mathrm{max}\{c :\:\exists s_2,\ldots,s_n\:\mathrm{with}\:
c<s_2<s_3<\ldots<s_n\\ \mathrm{and}\:[c,s_2,\ldots,s_n]\notin\mathfrak{q}_{s-m}
 \mathrm{\underline{and}}\: [c,s_2,\ldots,s_n]\: \mathrm{is \:
good} \}
\end{eqnarray*}
(then $1\le l_1< n+1$) and inductively,
\begin{eqnarray*}
L_i:=\mathrm{max}\{c :\:\exists s_{i+1}',\ldots,s_{n}'\:\mathrm{with}\:
L_1<\ldots<L_{i-1}<c<s_{i+1}'<\ldots<s_n'\\
\mathrm{and}\: [L_1,\ldots,L_{i-1},c,s_{i+1}',\ldots,s_n'] \: \mathrm{is
\: good} \\
\mathrm{and}\:[L_1,\ldots,L_{i-1},c,s_{i+1}',\ldots,s_n']\notin\mathfrak{q}_{s-m}\}
.\end{eqnarray*} 
Then  $[L_1,\ldots,L_n]\notin \mathfrak{q}_{s-m}$ and is good (furthermore
$L_n >n$ since $[1,\ldots,n] \in \mathfrak{q}_{s-m}$). I can write
$[L_1,\ldots,L_1]=[L^{\alpha}_1,\ldots,L^{\alpha}_h;L^{\beta}_1,\ldots,L^{\beta}_{n-h}]$ 
with\\ $\{ L^{\alpha}_1,\ldots,L^{\alpha}_h\} \cap \{
L^{\beta}_1,\ldots,L^{\beta}_{n-h} \} =\emptyset$. Choose
$b\in(\bigcap\limits_{i=s-m+1}^{s}\mathfrak{q}_i)\backslash\mathfrak{q}_{s-m}$ and
perform a base change in $A^{2n}$ (preserving the symmetry) by adding $b$ times the
$L_i^{\alpha}$ column of $\alpha$ to the $L_i^{\alpha}$ column of $\beta$, for
$i=1,\ldots,h$. Then
\[
[n+1,\ldots,2n]_{\mathrm{new}}=[n+1,\ldots,2n]_{\mathrm{old}}\pm
b^{h}[L_1,\ldots,L_n]_{\mathrm{old}}+b\mu,
\]
where $\mu\in\mathfrak{q}_{s-m}$ by the defining maximality property of the
$L$'s. Thus $[n+1,\ldots,2n]_{\mathrm{new}}\notin\mathfrak{q}_i$ for
$i=s-m,\ldots,s$, which is the inductive step. Therefore
after a sequence of base changes that preserve the symmetry
$\alpha\beta^t=\beta\alpha^t$ (and leave $\det(\alpha)$ unaltered) $\det(\beta)$ can
be made an
$A/(\det(\alpha))-$regular element, i.e. $\det(\alpha),\:\det(\beta)$ is an
$A-$regular sequence, which proves the lemma.      
\end{proof}

\vspace{3cm}
\begin{center} \bfseries References\end{center}
\mdseries
\small
\begin{itemize}
\item[\textrm{[A-C-G-H]}] E. Arbarello, M. Cornalba, P. Griffiths, J. Harris:
\emph{Geometry of Algebraic Curves}, Springer-Verlag, New York, N.Y. (1985) 

\item[\textrm{[Beau]}] A. Beauville: \emph{Complex algebraic surfaces}, London
Mathematical Society Student Texts 34, CUP (1996)

\item[\textrm{[Bom]}] E. Bombieri: \emph{Canonical models of surfaces of general type},
I.H.E.S. Publ. Math. 42 (1973), p. 171-219

\item[\textrm{[B-He]}] W. Bruns, J. Herzog: \emph{Cohen-Macaulay rings}, CUP (1998)

\item[\textrm{[B-V]}] W. Bruns, U. Vetter: \emph{Determinantal rings}, Springer LNM
1327 (1988)

\item[\textrm{[Cat1]}] F. Catanese: \emph{Babbage's conjecture, contact of
surfaces, symmetric determinantal varieties and applications}, Inv.
Math. 63 (1981), pp.433-465

\item[\textrm{[Cat1b]}] F. Catanese: \emph{On the moduli spaces of surfaces of general
type}, J. Differential Geom., 19 (1984), p. 483-515
 
\item[\textrm{[Cat2]}] F. Catanese: \emph{Commutative
algebra methods and equations of regular surfaces}, Algebraic
Geometry-Bucharest 1982, Springer LNM 1056 (1984), pp.68-111

\item[\textrm{[Cat3]}] F.Catanese: \emph{Equations of pluriregular varieties of general
type}, Geometry today-Roma 1984, Progr. in Math. 60, Birkh\"auser
(1985), pp. 47-67)

\item[\textrm{[Cat4]}] F. Catanese: \emph{Homological Algebra and Algebraic
surfaces}, Proc. Symp. in Pure Math., Volume 62.1, 1997, pp. 3-56

\item[\textrm{[Cil]}] C. Ciliberto: \emph{Canonical
surfaces with
$p_g=p_a=5$ and
$K^2=10$}, Ann. Sc.  Norm. Sup. s.IV, v.IX, 2 (1982), pp.287-336

\item[\textrm{[C-E-P]}] C. de Concini, D. Eisenbud, C. Procesi: \emph{Hodge algebras},
Ast\'{e}risque 91 (1982)

\item[\textrm{[C-S]}] C. de Concini, E. Strickland: \emph{On the variety of complexes},
Adv. in Math. 41 (1981), pp. 57-77 

\item[\textrm{[D-E-S]}] W. Decker, L. Ein, F.-O. Schreyer: \emph{Construction of
surfaces in $\mathbb{P}^4$}, J. Alg. Geom. 2 (1993), pp.185-237

\item[\textrm{[Ei]}] D. Eisenbud:  \emph{Commutative Algebra with a view towards
Algebraic Geometry}, Springer G.T.M. 150, New York (1995)

\item[\textrm{[E-P-W]}] D. Eisenbud, S. Popescu, C. Walter:  \emph{Symmetric locally free resolutions of coherent sheaves}, in preparation

\item[\textrm{[E-U]}] D. Eisenbud, B. Ulrich: \emph{Modules that are Finite
Birational Algebras}, Ill. Jour. Math. 41, 1 (1997), pp. 10-15

\item[\textrm{[En]}] F. Enriques: \emph{Le superficie algebriche},
Zanichelli, Bologna, 1949

\item[\textrm{[Gra]}] M. Grassi:
\emph{Koszul modules and Gorenstein algebras}, J. Alg. 180 (1996),
pp.918-953

\item[\textrm{[Gri-Ha]}] P. Griffiths, J. Harris: \emph{Principles of Algebraic
Geometry}, Wiley, New York (1978)

\item[\textrm{[Groth]}] A. Grothendieck: \emph{Local Cohomology}, Springer LNM 41
(1967)

\item[\textrm{[Ha]}] J. Harris: \emph{Algebraic Geometry, A First Course}, Springer
G.T.M. 133, New York (1993)

\item[\textrm{[Har]}] R. Hartshorne: \emph{Algebraic Geometry}, Springer G.T.M. 52, New
York (1977)

\item[\textrm{[Horr]}] G. Horrocks: \emph{Vector bundles on the
punctured spectrum of a local ring}, Proc. London Math. Soc. (3) 14
 (1964), pp. 689-713

\item[\textrm{[dJ-vS]}] T. de Jong, D. van Straten: \emph{Deformation of
the normalization of hypersurfaces}, Math. Ann. 288 (1990), pp.527-547

\item[\textrm{[Lip]}] J. Lipman: \emph{Dualizing sheaves, differentials and residues on
algebraic varieties}, Ast\'{e}risques 117 (1984)

\item[\textrm{[M-P]}] D. Mond, R. Pellikaan: \emph{Fitting
ideals and multiple points of analytic mappings}, Springer LNM 1414
(1987), pp.107-161

\item[\textrm{[P-S1]}] C. Peskine, L. Szpiro: \emph{Dimension projective finie et
cohomologie locale}, Publ. Math. IHES 42 (1972)

\item[\textrm{[P-S2]}] C. Peskine, L. Szpiro: \emph{Liaison des vari\'{e}t\'{e}s
alg\'{e}briques}, Inv. Math. 26 (1972), pp. 271-302

\item[\textrm{[Ro\ss ]}] D. Ro\ss berg: \emph{Kanonische Fl\"achen mit $p_g=5,\:q=0$
und
$K^2=11,12$}, Doktorarbeit Bayreuth, November 1996

\item[\textrm{[Sern]}] E. Sernesi: \emph{L'unirazionalit\`a della variet\`a dei 
moduli delle curve di genere dodici}, Ann. Scuola Norm. Pisa 8 (1981), pp.
405-439

\end{itemize}

\vspace{0.5cm}
\flushleft{\emph{Author's address:}}\\
\textsc{Mathematisches Institut der Universit\"at Bayreuth}\\
\textsc{Universit\"atsstra\ss e, D-95440 Bayreuth, Germany}\\
\vspace{0.2cm}
\emph{E-mail:} \ttfamily boehning@btm8x5.mat.uni-bayreuth.de 
\rmfamily 
    
\end{document}